\documentclass{amsart}
\usepackage{amssymb,amsxtra}
\usepackage{graphicx}
\usepackage{amscd}
\usepackage{amsmath}
\usepackage{amsfonts}
\usepackage{amssymb}
\newcommand{\lb}{\langle}
\newcommand{\rb}{\rangle}
\newcommand{\io}{\iota}
\newcommand{\al}{\alpha}
\newcommand{\C}{\mathbb C}

\title[fell bundles associated to groupoid morphisms]{Fell bundles associated to groupoid morphisms}
\author{Valentin Deaconu}
\address{Department of Mathematics\\ University of Nevada\\ Reno NV
89557-0084, USA}
\email[Valentin Deaconu]{vdeaconu@unr.edu}
\author{Alex Kumjian}
\email[Alex Kumjian]{alex@unr.edu}
\author{Birant Ramazan}
\email[Birant Ramazan]{ramazan@unr.edu}

\date{\today}
\begin{document}

\begin{abstract}
 Given a continuous open surjective morphism $\pi :G\rightarrow H$ of 
\'etale groupoids with amenable kernel, we construct a Fell bundle $E$
over $H$ and prove that its 
C*-algebra $C^*_r(E)$ is isomorphic to $C^*_r(G)$. This is related to
results of Fell concerning C*-algebraic bundles over groups. The case
$H=X$, a locally compact space, was treated earlier by Ramazan.  We
conclude that  
$C^*_r(G)$ is strongly Morita equivalent to a crossed product, the
C*-algebra of a Fell bundle arising from an action of the groupoid $H$
on a C*-bundle over $H^0$. We apply the theory to groupoid morphisms
obtained from extensions of dynamical systems and from morphisms of
 directed graphs with the path lifting property.
We also prove a structure theorem for abelian Fell bundles.
\end{abstract}

\maketitle

\S 1. INTRODUCTION

\bigskip

In his Memoir \cite{F1}, J. M. G. Fell generalizes Mackey's theory of
unitary representations of group extensions to a natural enrichment of
the concept of Banach *-algebra, called Banach *-algebraic
bundle. Given a normal subgroup $K$ of $G$, he constructs a bundle $B$
over $H = G/K$ with the fiber over the neutral element identified with
the algebra $L^1(K)$, such that $L^1(G)$ is isomorphic to the
cross-sectional algebra of $B$.  
He also proves that there is a one-to-one 
correspondence between isomorphism classes of Banach *-algebraic
bundles with one-dimensional fibers over the group $H$  and the family
of isomorphism classes of central topological extensions of $H$ by the unit
circle ${\mathbb T}$. 

Replacing $L^1(K)$ by $C^*(K)$, we get the notion of C*-algebraic
bundle over a locally compact group (see \cite[\S 11]{F2}). This may
be thought of as a continuous version of a group grading  
in a C*-algebra; one may regard the associated C*-algebra as a fairly 
general sort of crossed product of the fiber algebra over the neutral 
element by the group (in \cite{LPRS} it is shown that the C*-algebra is
endowed with a coaction by the group).  

There is a natural generalization
of the notion of C*-algebraic bundle to groupoids (see \cite{Y}),
which when specialized 
to topological spaces yields the more usual 
notion of (continuous) C*-bundle.  Such objects are often called 
Fell bundles (see \cite{K2,M}). 

The main result of this note is the construction of a Fell bundle
associated to certain groupoid homomorphisms.  We restrict ourselves
to \'etale groupoids (the range and the source maps are local
homeomorphisms), but  we expect that our results hold more generally. 
Given a continuous open surjective morphism $\pi : G\rightarrow H$ of \'etale
 groupoids with amenable kernel, we construct a Fell bundle $E$ over
 $H$ by extending a result of Ramazan's dissertation (see \cite{Ram}),
which appears in \cite{LR}.  The authors show that if $H$ is a 
 locally compact space, there is a C*-bundle over $H$ with fibers
given by C*-algebras associated to the fibers of $\pi$.  It follows
easily that $C^*_r(G)$ is isomorphic to the C*-algebra of continuous
sections of this bundle.   In our situation, their result may be applied
to the restriction of $\pi$ to the preimage of $H^0$ (which is an
\'etale groupoid) to obtain a C*-bundle over $H^0$.  This forms the
``nucleus'' of the desired Fell bundle $E$ (that is, its restriction to $H^0$).

Lee showed in \cite[Theorem 4]{L} that if $A$ is a C*-algebra and $X$
is a locally 
compact space, then $A$ may be realized as the C*-algebra associated to a
C*-bundle over $X$ if and only if there is a continuous open surjection
$\pi : \text{Prim}\,A \to X$.  In this case the fiber of the bundle over $x
\in X$ is the quotient of $A$ corresponding to the closed set
$\pi^{-1}(x) \subset  \text{Prim}\,A$.
Hence, if $\pi : Y \to X$ is a continuous open surjection of locally compact spaces, then 
$C_0(Y)$ is realizable as the C*-algebra associated to a
C*-bundle over $X$ with fibers $C_0(\pi^{-1}(x))$.
In Ramazan's result, the space $Y$ is replaced by a
groupoid and $\pi$ is required to be a groupoid morphism (where $X$ is
regarded as a  groupoid). In our result, $X$ is also replaced
by a groupoid, but we need Fell bundles rather than C*-bundles.

Several examples are considered, coming from extensions of dynamical systems and graph morphisms.

In the last section we consider abelian Fell bundles (the fibers over
the unit space are abelian C*-algebras) and prove a structure theorem
that states that every such bundle arises from a twisted groupoid covering.
We also give various examples of abelian Fell bundles and an
application of the structure theorem.


\bigskip

\S 2. FELL BUNDLES OVER GROUPOIDS

\bigskip
Recall the definition of a Fell bundle over a groupoid $G$ 
(see \cite{K2}).  Note that 
it follows from the first nine axioms that $E_u$ is a C*-algebra for
all $u \in G^0$, so it makes sense in (10) to require the positivity
of $e^*e$ for all $e \in E$. 

\noindent{\bf 2.1 Definition}. Let $G$ be a locally compact Hausdorff groupoid with
unit space $G^0$, range and source maps $r,s$ and set of composable
pairs $G^2$, which admits a left Haar system. A Banach bundle
$p:E\rightarrow G$ is said to be a {\em Fell bundle} if there is a continuous
multiplication $E^2\rightarrow E$, where 
\[E^2=\{(e_1,e_2)\in E\times E\;\mid \; (p(e_1),p(e_2))\in G^2\},\]and
an involution $e\mapsto e^*$ which satisfy the following axioms ($E_g$ is the fiber $p^{-1}(g)$).
\begin{itemize}

\item[1.] $p(e_1e_2)=p(e_1)p(e_2) \;\; \forall\;(e_1,e_2)\in E^2$;

\item[2.] the induced map $E_{g_1}\times E_{g_2}\rightarrow E_{g_1g_2}, \;(e_1,e_2)\mapsto e_1e_2$ is bilinear $\;\forall\; (g_1,g_2)\in G^2$;

\item[3.] $(e_1e_2)e_3=e_1(e_2e_3)$ whenever the multiplication is defined;

\item[4.] $||e_1e_2||\le ||e_1||\; ||e_2||\; \;\forall \;(e_1,e_2)\in E^2$;

\item[5.] $p(e^*)=p(e)^{-1}\;\; \forall \;e\in E$;

\item[6.] the induced map $E_g\rightarrow E_{g^{-1}}, \; e\mapsto e^*$ is conjugate linear for all $g\in G$;

\item[7.] $e^{**}=e\;\; \forall \;e\in E$;

\item[8.] $(e_1e_2)^*=e_2^*e_1^*\;\;  \forall \;(e_1,e_2)\in E^2$;

\item[9.] $||e^*e||=||e||^2\;\;  \forall \;e\in E$;

\item[10.] $e^*e\ge 0\;\;  \forall \;e\in E$.

\end{itemize}
\noindent A Fell bundle $E$ is called saturated if $E_{g_1}\cdot E_{g_2}$ is
total in $E_{g_1g_2}$ for all $(g_1,g_2)\in G^2$.

\noindent{\bf 2.2 Facts}.  For $g\in G$,  $E_{s(g)},\; E_{r(g)}$ are C*-algebras, and $E_g$
is a right Hilbert $E_{s(g)}$-module with inner product $\langle e_1,
e_2\rangle_s=e_1^*e_2$ and a left Hilbert $E_{r(g)}$-module with inner
product $\langle e_1, e_2\rangle_r=e_1e_2^*$. If $E$ is saturated, then $E_g$
is an $E_{r(g)}$-$E_{s(g)}$ equivalence bimodule, and for all $(g_1,g_2)\in G^2$, 
multiplication induces an isomorphism
$E_{g_1}\otimes_{E_u}E_{g_2}\cong E_{g_1g_2}$, where
$u=s(g_1)=r(g_2)$.   

The restriction
$E^0=E|_{G^0}$ is a C*-algebra bundle, and $C_0(E^0)$, the set of
continuous sections vanishing at $\infty$,  is a C*-algebra. 
We refer to Addendum 2 in \cite{K1} for other facts about C*-bundles.

Recall that a subset $S$ of a groupoid $G$ is called a  {\em bisection} if the
restrictions of the range and source maps to $S$ are injective.

\noindent{\bf 2.3 Lemma}. Let $E$ be a saturated Fell bundle over $G$ and let $U$
be an open bisection of $G$. Then the completion of $C_c(U,E)$, the continuous
compactly supported sections on $U$, with respect to the the supremum
norm, is an $A-B$ equivalence bimodule, where $A=C_0(r(U),E)$ and
$B=C_0(s(U),E)$, when endowed with the natural inner products and
actions. 

{\em Proof}.  The right and left multiplications are given by
\[(\xi\cdot b)(g)=\xi(g)b(s(g)),\;\; (a\cdot\xi)(g)=a(r(g))\xi(g),\]for $a\in A, b\in B, \xi\in C_c(U,E), g\in U.$
The inner products are
\[\langle \xi,\eta\rangle_B(s(g))=\xi(g)^*\eta(g), \;\;\langle \xi,\eta\rangle_A(r(g))=\xi(g)\eta(g)^*,\]
where $\xi,\eta\in C_c(U,E),\ g\in U$.
The positivity of the inner products follows from Definition 2.1.10.
$\Box$

For $G$ an \'etale groupoid and $p:E\rightarrow G$ a Fell bundle, one
can define multiplication and involution on the space of
compactly supported continuous sections $C_c(E)$ by 
\[(\xi\eta)(g)=\sum_{g=g_1g_2}\xi(g_1)\eta(g_2), \quad
\xi^*(g)=\xi(g^{-1})^*\] 
for $\xi, \eta \in C_c(E)$.  
Define an inner product  $\langle \xi, \eta\rangle=P(\xi^*\eta)$
for $\xi, \eta \in C_c(E)$, where $P:C_c(E)\rightarrow C_c(E^0)$ is
the restriction map;  denote by $L^2(E)$ the completion of $C_c(E)$ in
the norm defined by this inner product (so  
$\|\xi\|^2 = \|\langle \xi, \xi\rangle\|$).
Observe that $C_c(E)$ acts by left multiplication on $L^2(E)$.

\noindent{\bf 2.4 Definition}. The C*-algebra $C^*_r(E)$ is defined as
the completion of $C_c(E)$  in 
${\mathcal L}(L^2(E))$, with respect to the operator norm.

\noindent{\bf 2.5 Remark}. The restriction map
$C_c(E)\rightarrow C_c(E^0)$ extends to a faithful conditional expectation
$P:C^*_r(E)\rightarrow C_0(E^0)$.  
Also, $L^2(E)$ is the Hilbert module associated to a
bundle of Hilbert modules $V$ over $G^0$, where $\displaystyle
V_u=\bigoplus_{s(g)=u}E_g$ (see \cite[3.3]{K2}). If $E$ is saturated,
there is natural action of the groupoid $G$ on the C*-algebra
bundle ${\mathcal K}(V)$ with fibers 
${\mathcal K}(V_u)\cong V_u\otimes V_u^*$.
One can form the semi-direct product bundle 
$G\ltimes {\mathcal K}(V)$ over $G$.  Kumjian proved in \cite[4.5]{K2}
that, if $E$ is saturated, then
$C^*_r(G\ltimes{\mathcal K}(V))$ and $C^*_r(E)$ are strongly Morita
equivalent.

 \bigskip

 \S 3. THE FELL BUNDLE ASSOCIATED TO A GROUPOID MORPHISM

\bigskip

\noindent{\bf 3.1 Definition}. 
Let $G$ and $H$ be topological groupoids.
A {\em groupoid morphism} $\pi :G\rightarrow H$ is a
continuous map that
intertwines both the range and source maps and that satisfies
$\pi(g_1g_2)=\pi(g_1)\pi(g_2)$ for all $(g_1,g_2)\in G^2$.  It follows that $\pi(G^0)\subset H^0$.

\noindent{\bf 3.2 Definition}. A {\em groupoid fibration} is an
 open surjective morphism 
of locally compact groupoids $\pi :G\rightarrow H$ with the property
that for any $h\in H$ and $x\in G^0$ with $\pi(x)=s(h)$ there is $g\in
G$ with $s(g)=x$ and $\pi(g)=h$. 
If $g$ is unique for any such $h$ and $x$, then $\pi$ is called 
 a {\em groupoid covering}. Note that for a groupoid covering we have
 $\pi^{-1}(H^0)=G^0$ (see \cite{Br}).

For example, if $Y\rightarrow X$ is a Serre
fibration of topological spaces, then $\pi_1(Y)\rightarrow \pi_1(X)$ is a groupoid
fibration, and if $\tilde{X}\rightarrow X$ is a covering, then
$\pi_1(\tilde{X})\rightarrow \pi_1(X)$ is a groupoid covering. Here
$\pi_1(X)$ denotes the fundamental groupoid of the space $X$. 

The following is essentially a restatement of Ramazan's result.

\noindent{\bf 3.3 Lemma}. Given an open surjective morphism $\pi:K\rightarrow X$,
where $K$ is a locally compact amenable groupoid and $X$ is a locally
compact space, there is a C*-bundle $F$ over $X$ with fibers
$C^*_r(K(x))$, where $K(x)=\pi^{-1}(x)$. Moreover, $C^*_r(K)$ is
isomorphic to the C*-algebra of continuous sections of this bundle. 

{\em Proof}. We take $F$ to be the disjoint union of $C^*_r(K(x))$
over $x\in X$ with the bundle structure defined as in Proposition 1.6
in \cite{F1} with $\Gamma = C_c(K)$; 
we may view $\Gamma$ as sections of this bundle by means of the canonical maps $\pi_x:  C_c(K) \to C^*_r(K(x))$. The continuity of the norm is
proved in Th\'eor\`eme 
2.4.6 in Ramazan's thesis (\cite{Ram}) or Corollary 5.6 in 
\cite{LR}. Using Proposition 1.7 in \cite{F1} with  $\Gamma$ the compactly supported sections of $F$, we get that $C^*_r(K)\cong
C_0(F)$.$\;\Box$

\noindent{\bf 3.4 Theorem}. Given an open surjective morphism of \'etale  groupoids
$\pi:G\rightarrow H$ with amenable kernel $K := \pi^{-1}(H^0)$, there is
 a Fell bundle $E=E(\pi)$ over $H$ such that
$C^*_r(G)\cong C^*_r(E)$. Moreover, $C^*_r(G)$ is strongly Morita
equivalent to a crossed product $C^*_r(H\ltimes {\mathcal K}(V))$ 
(see  Remark 2.5).

{\em Proof}. 
Using the lemma for the restriction of $\pi$ to $K=\ker \pi=\pi^{-1}(H^0)$, we get a C*-bundle $F$ with fibers $C^*_r(K(x))$ over the unit space $H^0$. We will  extend 
this C*-algebra bundle to a Fell bundle $E=E(\pi)$ over  $H$.
 
Note that $K$ is an open  \'etale subgroupoid of $G$,
and there is a faithful conditional expectation 
$\Phi: C^*_r(G)\rightarrow C^*_r(K)$. Indeed, if $\Psi: C^*_r(K)\to
C_0(K^0)$ is the canonical conditional expectation, where $K^0=G^0$ is
the unit space, then $\Phi$ must be faithful since $\Psi\circ\Phi$ is
faithful. 
We construct a Hilbert module $M(\Phi)$ over
$C^*_r(K)$ by completing $C_c(G)$ with respect to the norm given by
the inner product $\langle f_1,f_2\rangle=\Phi(f_1^*f_2)$ for
$f_1,f_2\in C_c(G)$. The right multiplication is given by convolution
with elements in $C_c(K)\subset C^*_r(K)$ (and extending by
continuity). Since $C^*_r(K)\cong C_0(F)$ and $F$ is fibered over
$H^0$, it follows that $M(\Phi)\cong C_0(B)$, where $B$ is a bundle of
Hilbert modules over $H^0$, with  $B_x$ a Hilbert module   over
$F_x=C^*_r(K(x))$ for each $x\in H^0$ (see 1.7 in \cite{K2}). The
fiber $B_x$ is the completion of 
$C_c(G(x))$, where $G(x)=\{g\in G \mid s(\pi(g))=x\}$ and the inner product is the natural restriction of the above inner product on $C_c(G)$.
We define $E_x=F_x$ for $x\in H^0$, and for arbitrary $h\in H$, we define $E_h$ to be the
completion of $C_c(\pi^{-1}(h))$ in the norm coming from the inclusion $C_c(\pi^{-1}(h))\subset C_c(G(s(h)))\subset B_{s(h)}$.  Note that $E_h$ is a submodule of $B_{s(h)}$. The 
 multiplication $E_{h_1}\times E_{h_2}\rightarrow E_{h_1h_2}$ is
 defined by  
 \[(\xi\eta)(g)=\sum_{g_1g_2=g}\xi(g_1)\eta(g_2), \quad\text{for}\quad\xi\in
 C_c(\pi^{-1}(h_1)), \; \eta\in C_c(\pi^{-1}(h_2));\] 
 for $\xi\in C_c(\pi^{-1}(h))$ we define $\xi^*\in C_c(\pi^{-1}(h^{-1}))$ by
 $\xi^*(g)=\overline{\xi(g^{-1})}$.  Observe that the norm on $E_h$ inherited from $B_{s(h)}$ satisfies 
 $||\xi||=||\xi^*\xi||^{1/2}$  for $\xi\in E_h$. Moreover, the element $\xi^*\xi\in E_{s(h)}$ is positive since $\xi^*\xi=\langle \xi,\xi\rangle_{s(h)}$, where $\langle\cdot,\cdot\rangle_{s(h)}$ is the inner product on $B_{s(h)}$.  The bundle structure for the union of $E_h$'s is given by $C_c(G)$, using the fact that each element
in $C_c(\pi^{-1}(h))$ is the restriction of an element in $C_c(G)$ (see Proposition 10.7 in \cite{F2}). To prove the continuity of the norm, fix $h_0\in H$ and take $U$ an open bisection of $H$ containing $h_0$. For $\xi\in C_c(G)$, denote by $\xi_h$ the restriction of $\xi$ to $\pi^{-1}(h)$. By a partition of unity argument, it will be sufficient to consider $\xi$ with support in $\pi^{-1}(U)$.  Since $||\xi_h||=||\xi_h^*\xi_h||^{1/2}$   and $\xi_h^*\xi_h=(\xi^*\xi)_{s(h)}$, it follows that  the map $x\mapsto ||(\xi^*\xi)_x||$ on $s(U)$ is continuous.   It is straightforward to check all the
other axioms of a Fell bundle. The  bundle $E$ is always saturated.

We can now identify $C^*_r(K)$ with $C_0(E^0)$, where $E^0$
is the restriction of $E$ to $H^0$. To prove that  $C^*_r(G)\cong C^*_r(E)$, we use the natural extension of the
map $\psi: C_c(G)\rightarrow C_c(E)$ given by 
$$\psi(f)(h)=f\mid_{\pi^{-1}(h)} \in C_c(\pi^{-1}(h)) \subset E_h$$ 
to get an isomorphism $U_{\psi}$ between the $C_0(E^0)$-Hilbert
modules $M(\Phi)$ and  $L^2(E)$ (see \S 2 before Definition 2.4). Indeed, the module structures and the
inner product are preserved since in each case both are derived from
convolution and involution on $C_c(G)$ (note that $\psi$ is a map of $*$-algebras).

Both $C^*_r(G)$ and $C^*_r(E)$ are represented on these isomorphic Hilbert modules, using the left regular representation.
The same map $\psi$ preserves the product, and it induces an
isomorphism $\alpha_{\psi}= \text{Ad}\, U_{\psi}$ between these C*-algebras.

The last part of the statement follows from Kumjian's result
mentioned in Remark 2.5.$\;\Box$

\noindent{\bf 3.5 Example}. Consider $G$ a discrete group, $K$ a normal  subgroup and let $H=G/K$
with $\pi :G\to H$ the canonical morphism. Then we get a Fell bundle $E$ over $H$ with
the fiber $C_r^*(K)$ over the identity element, such that $C_r^*(G)\cong C^*_r(E)$.  This is a particular case of the construction done by Fell in the context of homogeneous Banach *-algebraic bundles over locally compact groups (see Example 3 on page 77 in \cite{F1}). Recall that a Fell bundle over a discrete group is equivalent to a grading.

We specialize to the discrete $3$-dimensional Heisenberg group $G\subset SL(3,{\mathbb Z})$. The group $G$ consists  of matrices of the form 
\[[a,b,c]:=\left[\begin{array}{ccc}1&a&c\\0&1&b\\0&0&1\end{array}\right].\] 
The group operation is
\[[a,b,c][a',b',c']=[a+a',b+b',c+c'+ab'].\]
It is easy to show that $G$ is an extension of ${\mathbb Z}^2$ by its
center ${\mathbb Z}$.  We have a surjective homomorphism
$\pi:G\to{\mathbb Z}^2$ given by $\pi[a,b,c]=(a,b)$ with
$\ker\pi\cong{\mathbb Z}$. The C*-algebra $C^*(G)\cong C^*_r(G)$ is
called the rotational algebra in \cite {AP}, and it may be understood
as the algebra of continuous sections of a field of C*-algebras over
the unit circle ${\mathbb T}$. Our construction from the morphism
$\pi$ gives a new perspective: $C^*(G)$ is the C*-algebra of a Fell
bundle over the group ${\mathbb Z}^2$ with fibers isometric to $C({\mathbb
  T})$.  Note that the complexity of the structure of $C^*(G)$ is
contained in the definition of the product between the fibers, since
for instance $C({\mathbb T}^3)\cong C^*({\mathbb Z}^3)$ is also  the
C*-algebra of a Fell bundle over ${\mathbb Z}^2$ with fibers
isometric to $C({\mathbb T})$. 
 
There is another characterization of $G$ as an extension of ${\mathbb
  Z}$ by ${\mathbb Z}^2$ coming from the semi-direct product
decomposition $G\cong {\mathbb Z}^2\rtimes{\mathbb Z}$. Here ${\mathbb
  Z}^2$ is generated by $[1,0,0]$ and $[0,0,1]$, and the action of
${\mathbb Z}$ is defined by conjugation with $[0,1,0]$. The morphism
$\pi: G\to{\mathbb Z}$ given by $\pi[a,b,c]=b$ describes $C^*(G)$ as
  the C*-algebra of  a Fell bundle over ${\mathbb Z}$ with fibers
  isometric to  $C({\mathbb T}^2)$.   

\noindent{\bf 3.6 Example}. For $G$ an \'etale  groupoid, consider the equivalence relation \[R=\{(x,y)\in G^0\times G^0 \mid  \exists \; g\in G \;\mbox{such that}\; r(g)=x, s(g)=y\}\] and 
the map $\pi:G\rightarrow R, \;\;\pi(g)=(r(g),s(g))$. Assume that the
isotropy group bundle is amenable and open in $G$. Then $\pi$ is an
open surjective morphism, and $C^*_r(G)$ may be realized as the C*-algebra
of a Fell bundle over $R$.

\bigskip


\S 4. FELL BUNDLES FROM GRAPH MORPHISMS

\bigskip

\noindent{\bf 4.1 Definition}. Let $V,W$ be (finite) graphs. A graph morphism $\phi:V\rightarrow W$ is a map which preserves the incidences. If $\phi$ is surjective, we say that
it has the path lifting property (or that it is a fibration)  if for any vertex $v\in V^0$ and any edge $b\in W^1$ starting at $w=\phi(v)$ there is an edge $a\in V^1$ starting at $v$  with $\phi(a)=b$.

Recall (see \cite{K3}) that, if $V$ has no sinks,
 the graph C*-algebra $C^*(V)$ is the C*-algebra of the amenable groupoid  
\[G_V=\{(a,p-q,a')\in X_V\times\mathbb{Z}\times X_V\mid\; \sigma^p(a)=\sigma^q(a')\},\]
where $\sigma(a_1a_2a_3\cdots)=a_2a_3a_4\cdots$ is the shift map and $X_V$ is the space  of infinite paths $a_1a_2a_3\cdots$ of concatenated edges in $V^1$.

\noindent{\bf 4.2 Proposition}. Assume that $V$ and $W$ have no sinks. A graph
morphism  $\phi$  with the path lifting property induces a continuous open surjection
\[\varphi: X_V\rightarrow X_W,\quad\text{given by}\quad\varphi(a_1a_2a_3
\cdots)=\phi(a_1)\phi(a_2)\phi(a_3)\cdots\] 
between the  infinite path
spaces, and  an open  surjective morphism
\[\pi: G_V\rightarrow G_W, \quad\text{given by}\quad
\pi(a,k,a')=(\varphi(a),k,\varphi(a'))\] 
between the associated groupoids, which is a fibration with kernel 
\[K=\{(a,0,a')\in G_V\;\mid\; \varphi(a)=\varphi(a')\}.\] 
Hence, $C^*_r(G_V)$ is isomorphic to the C*-algebra of a Fell bundle over $G_W$.

{\em Proof}.  To show that
$\varphi$ is surjective, consider an infinite path $b_1b_2\cdots \in X_W$ beginning
 at $w_1\in W^0$. Since $\phi$ is onto, there is $v_1\in V^0$ with $\phi(v_1)=w_1$. By the path lifting property, there is $a_1\in V^1$ with $\phi(a_1)=b_1$. Continuing inductively, it follows that there is $a_1a_2\cdots\in X_V$ such that 
 $\varphi(a_1a_2\cdots)=b_1b_2\cdots$, and therefore $\varphi$ is surjective.
 Consider a cylinder set $Z\subset X_V$ determined by a finite path $a_1\cdots a_n$, i.e.
 \[Z=\{a_1\cdots a_nx_1x_2\cdots \in X_V\;\mid x_1x_2\cdots\in X_V\}.\]  Again by the path lifting property, $\varphi(Z)$ is the cylinder set in $X_W$ determined by the finite path $\phi(a_1)\cdots\phi(a_n)$. It follows that $\varphi:X_V\to X_W$ is continuous and open. 
 
By definition, we have 
\[\pi((a,k,a')(a',l,a''))=\pi(a,k+l,a'')=(\varphi(a),k+l,\varphi(a''))=\pi(a,k,a')\pi(a',l,a''),\]
and $\pi$ intertwines the range and source maps, therefore is a groupoid morphism. Since $\varphi$ is surjective and maps cylinder sets onto cylinder sets, it follows that $\pi$ is surjective, continuous and open. To show that $\pi$ is a groupoid fibration, consider $h=(b,k,b')\in G_W$ and $a'\in G_V^0=X_V$ with $\varphi(a')=s(h)=b'$. Since $\varphi$ is onto and intertwines   the shift maps, we can find $g=(a,k,a')\in G_V$ with $\pi(g)=h$. Hence $\pi$ is a groupoid fibration. 
Now $(a,k,a')\in \ker\pi$ iff $\varphi(a)=\varphi(a')$ and $k=0$, and the last statement of the proposition follows from  Theorem 3.4. $\Box$

\noindent{\bf 4.3 Example}. Consider the  graphs $V, W$ with $V^0=\{v\},\
W^0=\{w\},\ V^1=\{a,b,c\},\ W^1=\{1,2\}$ and $\phi(a)=\phi(b)=1,
\phi(c)=2$. Then $\phi$ induces a continuous map $\varphi:\{a,b,c\}^{\mathbb
N}\rightarrow \{1,2\}^{\mathbb N}$ between the infinite path spaces,
and a  morphism $\pi$ between the Cuntz 
groupoids $G_V$ and $G_W$. Hence, the Cuntz algebra ${ O}_3$ is
isomorphic to the C*-algebra of a Fell bundle over $G_W$. Note that
 for $x\in X_W$, \[K(x)=\{(y,0,z)\in
G_V\;\mid\; \varphi(y)=\varphi(y)=x\}.\] The fibers
$C^*(K(x))$ of the Fell bundle are isomorphic to  $M_{2^n}$, where $n$
is the number of $1$'s in $x$. For $n=\infty$, $M_{2^{\infty}}$ is the
UHF-algebra of type $2^{\infty}$. It is interesting to note that these
fibers are different.

\noindent{\bf 4.4 Example}. For a graph $V$ with no sinks consider the
collapsing map $\phi$ onto the graph $Z$ with  
one vertex and one loop. In this case $X_Z=\{*\}$ and
$G_Z\cong{\mathbb Z}$. The morphism $\pi$ is given by $\pi(x,k,y)=k$
and it induces the canonical $\mathbb{Z}$-grading on the C*-algebra $C^*(V)$.

\noindent{\bf 4.5 Remark}. Notice that the open map $\varphi :X_V\to X_W$ from Proposition 4.2 could be interpreted as a factor map between the topological Markov shifts  $(X_V,\sigma)$ and $(X_W,\sigma)$. More generally, 
let $X, Y$ be locally compact spaces and let $\sigma:X\to X$, 
$\tau:Y\to Y$ be two
local homeomorphisms. Assume that there
is a continuous  surjection $\varphi: X\rightarrow Y$ such that
$\tau\circ\varphi=\varphi\circ\sigma$. In the language of dynamical
systems, $(X,\sigma)$ is an extension of $(Y,\tau)$, or $(Y,\tau)$ is
a factor of $(X,\sigma)$. If $\varphi$ is also open, it induces an
 open surjective groupoid morphism 
\[\pi:G(\sigma)\rightarrow G(\tau)\quad\text{given by}\quad
 \pi(x,k,y)=(\varphi(x),k,\varphi(y)),\]
where $G(\sigma)=\{(x,m-n,y)\in X\times{\mathbb Z}\times X\mid
\sigma^mx=\sigma^ny\}$ and $G(\tau)$ is defined in the same way. These
groupoids are amenable by Proposition 2.4 in \cite{Re2}. 
Moreover, $\pi$ is a groupoid fibration.

\noindent{\bf 4.6 Example}.
 Let $X=\{1,2\}^{\mathbb N}\times {\mathbb  T}, Y=\{1,2\}^{\mathbb
 N}$  and let $\sigma : X \to X$ be given by
\[\sigma(a_1a_2\cdots,z)=\left\{\begin{array}{l}(a_2a_3\cdots,z^2)\quad\text{if}\quad
 a_1=1\\(a_2a_3.\cdots,z^3)\quad\text{if}\quad a_1=2\end{array}\right..\]
Let $\tau : Y \to Y$ and $\varphi:X\to Y$ 
be given by $\tau(a_1a_2\cdots)=a_2a_3\cdots$ and 
$\varphi(a,z)=a$.

In this case the C*-bundle over the unit space $G(\tau)^0=\{1,2\}^{\mathbb N}$ has fibers 
over $a\in \{1,2\}^{\mathbb N}$ isomorphic to Bunce-Deddens algebras of type $2^n3^m$ where $n, m\in {\mathbb N}\cup\{\infty\}$ are the number of $1's$ and $2's$ in $a$, respectively.

\noindent{\bf 4.7 Remark}. The notion of graph morphism with the path lifting
property  can be generalized to continuous graphs (see \cite{Ka}). This gives a larger
class of examples.

\bigskip

\S 5. ABELIAN FELL BUNDLES
\bigskip

\noindent{\bf 5.1 Definition}. A  Fell bundle over a
groupoid $H$ is called {\em abelian} if $E_u$ is an abelian C*-algebra for all
$u\in H^0$.

\noindent{\bf 5.2 Example}. Consider $\pi:G\to H$  a groupoid covering (see Definition 3.2). Then
the corresponding Fell bundle $E(\pi)$ is abelian, since
$\pi^{-1}(H^0)=G^0$.  So in the notation of Theorem 3.4,
$E_u = C_0(\pi^{-1}(u))$ for all $u \in H^0$.

Coverings of groupoids are intimately related to groupoid
actions on spaces. We will prove that every covering comes from such
an action and vice versa. 
Other examples are related to some groupoid extensions. A Fell line bundle is an abelian Fell bundle (see Example 5.5). The main
result in this section is a structure theorem for abelian Fell bundles
that, loosely speaking, asserts that every such bundle arises from a 
``twisted'' covering.

Recall (see \cite{MRW1})  that a groupoid $G$ is said to act
(on the left) on a locally compact space $X$, if there are given 
a continuous, open surjection $\rho : X \rightarrow G^0$
and a continuous map
\[G * X \rightarrow X, \quad\text{write}\quad (g , x)\mapsto
g \cdot x,\]
where
$G * X = \{(g , x)\in G \times X \mid s(g) = \rho (x)\},$
that satisfy

i) $\rho (g \cdot x) =r (g), \,\,\forall \,\,(g , x) \in G * X,$

ii) $(g _1, x) \in G * X,\,\, (g_2, g_1)
\in G ^2$ implies $(g _2g _1, x),
(g _2, g _1\cdot x) \in G * X$ and
\[g _2\cdot(g _1\cdot x) = (g _2g _1)
\cdot x,\]

iii) $\rho (x)\cdot x = x, \,\,\,\forall \,\, x\in X.$

Note that the fibered product $G * X$ has a natural
structure of  groupoid, called the semi-direct product or
action groupoid and denoted  by $G \ltimes X$ (cf.\ \cite{AR}), where
\[(G \ltimes X)^2 = \{ ((g_2, x_2),(g _1, 
x_1)) \mid \,\,x_2 = g _1\cdot x_1\}\]
\[(g _2, g_1\cdot x_1)(g _1,  x_1) = 
(g _2g _1, x_1)\]
\[(g,x)^{-1} = (g ^{-1}, g \cdot x).\]
Here the source and range maps are
 \[s(g ,x) = (g^{-1} g , x)=(s(g),x)=
(\rho ( x), x),\quad
r(g ,x) = (g g^{-1} ,g\cdot x)= (r(g), g\cdot x)=(\rho(g\cdot  x),g\cdot x),\]
and the unit space may be identified with $X$.

For $G$ an \'etale groupoid, consider the groupoid morphism $\pi:
G \ltimes X\rightarrow G,\quad \pi(g,x)=g$.  Then $\ker\pi$ may be
identified with $X$, and $C^*_r(G \ltimes X)$ may  be regarded as the
C*-algebra of an abelian Fell bundle over $G$ since $\ker\pi=X$ is amenable and
$E(\pi)_u \cong C_0(\rho^{-1}(u))$ for all $u \in G^0$.

\noindent{\bf 5.3 Proposition}.
Let $G$ be an \'etale groupoid acting on the locally
compact space $X$. Then the morphism $\pi: G\ltimes X\to G,\quad
\pi(g,x)=g$ is a covering, and gives rise to an abelian Fell bundle
over $G$. 
Moreover, every covering of \'etale groupoids is of this form. More
precisely, given a covering map $\pi: G \to H$ of \'etale groupoids,
there is a space $X$ and an action of $H$ on $X$ such that 
$G\cong H\ltimes X$. 

{\em Proof}. 
Let $g \in G$ and let $x \in (G\ltimes X)^0 = X$ such that 
$\pi(x) = s(g)$; then $(g, x) \in G \ltimes X$ satisfies 
$\pi(g, x) = g$ and $s(g, x) = x$, and is the unique such element in 
$G\ltimes X$.  Hence, $\pi: G\ltimes X\to G$ is a covering.

Comversely, let $\pi: G \to H$ be a covering of \'etale groupoids.
Set $X=G^0$ and $\rho=\pi^0: X\rightarrow H^0$. The map 
$\pi*s:G\rightarrow H \ltimes X$  is a continuous bijection 
(bijectivity follows from the definition of covering). To prove that
$\pi*s$ is open, consider 
\[{\mathcal U}=\{U\subset G\;\mid\;U\;{\mbox{is an open}}\;
\mbox{bisection of}\;G\;{\mbox{such that}}\; \pi(U)\;{\mbox{is an open
    bisection of}}\;H \}.\]
Then\ ${\mathcal U}$ forms a basis for $G$, and if $U\in{\mathcal U}$,
    then $(\pi*s)(U)$ is open in $H\ltimes X$. Hence,
$\pi*s$ is open and thus an isomorphism of locally compact groupoids. 
$\Box$

See \cite[Theorem 1.8]{KS} for an analogous result when $H$ is a
group.  The authors give conditions under which the groupoid
C*-algebra is Morita equivalent to a crossed product of an abelian
C*-algebra by an action of $H$.

\noindent{\bf 5.4 Example}.  Let $\pi : G \to H$ be a  open surjective morphism where $G$
and $H$ are \'etale groupoids.  
Suppose that the restriction of $\pi$ to $G^0$ induces a
homeomorphism $G^0 \cong H^0$ and that $A = \pi^{-1}(H^0)$ is a sheaf
of abelian groups over $G^0$.  Then the resulting sequence:
$$
A \overset{\io}{\longrightarrow} G \overset{\pi}{\longrightarrow} H
$$
is called an abelian extension where $\io$ is the inclusion map.
Since $A$ is amenable, 
the main theorem in \S 3 applies and we get a Fell 
bundle $E = E(\pi)$ over $H$.  Moreover, $E_u = C^*(A_u)$  for all $u \in H^0$,
where $A_u$ is the fiber of $A$ over $u$.  Since $A$ is a sheaf of
abelian groups, $E_u$ is abelian  for all $u \in H^0$.
Hence, $E$ is an abelian Fell bundle. 

\noindent{\bf 5.5 Example}. Let $G$  be a proper ${\mathbb T}$-groupoid over $H$,
that is, $G$ is a groupoid endowed with the structure of 
a principal ${\mathbb T}$-bundle over $H$ compatible with the groupoid structure (see \cite{K1}). 
Form the associated line bundle:
$$E =G *_{\mathbb T}{\mathbb C} = (G \times {\mathbb C}) / {\mathbb T}$$ 
(where  $t(g, z) = (t\cdot g, t^{-1}z)$ ).
One defines multiplication and involution as follows:
\begin{align*}
  (g_1, z_1)(g_2, z_2) &= (g_1g_2, z_1z_2), \\
                (g, z)^* &= (g^{-1},{\bar z} ).
\end{align*}
One verifies that $E$ is a Fell bundle over $H$ with these operations
and that $E_u \cong {\mathbb C}$ for all $u \in H^0$. A bundle of this type is called a \emph{Fell line bundle}.  Any Fell
bundle for which $E_h$ is one dimensional for all $h \in H$ is of this type.
Its C*-algebra is isomorphic to the twisted groupoid C*-algebra of $G$ (see \cite{MW}).

\noindent{\bf 5.6 Theorem}. Given a saturated abelian Fell bundle $E$ over an
\'etale groupoid $H$, there is a groupoid $G$, a covering $\pi:G\to H$
and a one-dimensional Fell bundle $L$ over $G$ such that
$C^*_r(L)\cong C^*_r(E)$. 

{\em Proof}. Set $X=\widehat{C_0(E^0)}$ i.e. $C_0(E^0)=C_0(X)$. Since
$C_0(X)$ is the C*-algebra of a bundle over $H^0$,  we get a continuous open
surjection $\rho:X\to H^0$ (see \cite{L}).  For  $u\in H^0$ we have an
isomorphism $j_u : E_u \cong C_0(X_u)$ where $X_u=\rho^{-1}(u)$;
for each $h \in H$, $E_h$ is a $C_0(X_{r(h)})-C_0(X_{s(h)})$ 
equivalence bimodule and hence we
get a homeomorphism $\alpha_h: X_{s(h)}\to X_{r(h)}$ (see
\cite[Corollary 3.33]{RW}, \cite[Corollary 6.27]{Ri}).  
This defines a map $H*X\to X, (h,x)\mapsto \alpha_h(x)$. We wish to
show that this defines an action of $H$ on $X$ (see the definition
following  Example 5.2). Conditions i) and iii) are immediate, while
condition ii) follows from the isomorphism $E_{h_1h_2}\cong
E_{h_1}\otimes E_{h_2}$.

By Lemma 2.3, if $U$ is an open bisection of $H$, the completion of
$C_c(U,E)$ may be endowed with the structure of an
 $A-B$ equivalence bimodule, where $A=C_0(r(U),E)$ and
$B=C_0(s(U),E)$, with the natural inner products and
actions.    Note that 
 $A \cong C_0(\rho^{-1}(r(U)))$ and
$B \cong C_0(\rho^{-1}(s(U)))$.
So again by \cite[Corollary 3.33]{RW}) there is a homeomorphism 
$\rho^{-1}(s(U)) \cong \rho^{-1}(r(U))$ compatible with the above
fiberwise homeomorphisms $\alpha_h$.  This proves that the map
$H*X\to X$ is continuous.
Thus, $H$ acts on $X$.  We set $G = H \ltimes X$ and  
by the above proposition, the map $\pi : H \ltimes X \to H$ given by 
$\pi(h, x) = h$ is a covering.

We construct $L$ piecewise as follows. 
Let  $h \in H$; as noted above,
$E_h$ is a $C_0(X_{r(h)})-C_0(X_{s(h)})$ equivalence bimodule
and the equivalence induces a homeomorphism 
$\al_h: X_{s(h)} \to X_{r(h)}$.  Then by Proposition A3 in \cite{Rae}
there is a Hermitian line bundle $L(h)$ over 
$\pi^{-1}(h) = \{ (h, x) \mid \rho(x) = s(h) \}$ and an isomorphism 
$\io_h : E_h \cong C_0(\pi^{-1}(h), L(h))$ such that for all 
$\xi, \eta \in E_h$ and $g \in \pi^{-1}(h)$ we have
\begin{align*}
j_{r(h)}(\xi\eta^*)(r(g)) &= \lb \io_h(\eta)(g), \io_h(\xi)(g) \rb  \\
j_{s(h)}(\xi^*\eta)(s(g)) &= \lb \io_h(\xi)(g), \io_h(\eta)(g) \rb,
\end{align*}
where $\lb \cdot, \cdot \rb$ denotes the fiberwise sesquilinear
inner-product (conjugate linear in the first variable).
For $(h_1, h_2) \in H^2$ with $u = s(h_1) = r(h_2)$, we have an
isomorphism 
$m : L(h_1) \otimes_{C_0(X_u)}  L(h_2) \cong  L(h_1h_2)$ such that
for $\xi_1 \in E_{h_1}$ and $\xi_2 \in E_{h_2}$ we have
$$
m(\io_{h_1}(\xi_1) \otimes\io_{h_2}(\xi_2)) = \io_{h_1h_2}(\xi_1\xi_2).
$$
Moreover, involution on $E$ defines a conjugate linear isomorphism
$L(h) \cong L(h^*)$ for all $h \in H$.  
Now for $g \in G$, we define  $L_g$  as the fiber over $g$ of the Hermitian line bundle $L(\pi(g))$ (note 
$L_g \cong \C$).  We use the above
operations to define multiplication and involution.

We wish to endow $L$ with the structure of a complex line bundle over $G$.
Let $g\in G$, and let $U$ be an open bisection of $H$ containing
$\pi(g)$. Then $C_0(\rho^{-1}(r(U)))$ and 
$C_0(\rho^{-1}(s(U)))$ are abelian and Morita equivalent
(and hence isomorphic as noted above).
Again by \cite[Proposition A3]{Rae} the equivalence bimodule may be
identified with $C_0(\pi^{-1}(U), L^U)$ where $L^U$ is a
Hermitian line bundle over $\pi^{-1}(U)$.  
One checks that $L(h) \cong L^U\big|_{\pi^{-1}(h)}$
for $h \in U$. 
This gives us the topology of $L$ as a line bundle over $G$. 
The product and involution defined above are compatible with this
topology and, hence, $L$ is a Fell bundle over $G$. 
It follows that $ C^*_r(L)\cong C^*_r(E)$, where the map 
$\psi: C_c(L) \to C_c(E)$ is given by 
$\psi(f)(h) = \io_h^{-1}\big(f|_{\pi^{-1}(h)}\big)$. $\Box$

Note that $E$ may be regarded as a push forward of $L$ under $\pi$. 
Specializing Example 5.4 to the case of groups,  the above theorem yields a result that is no doubt
well known to specialists.

\noindent{\bf 5.7 Corollary}. Let $\pi : G \to H$ be a surjective homomorphism, where $G$
and $H$ are discrete groups, such that $A=\ker \pi$ is abelian. Then we
get an action of $H$ on the dual $\hat{A}$ and a Fell line bundle over the
groupoid $H\ltimes\hat{A}$ defined by a two-cocycle
$\omega:(H\ltimes\hat{A})^2\to {\mathbb T}$ 
such that $C^*_r(G)\cong C^*_r(H\ltimes\hat{A},\omega)$. 

Note that the line bundle over $H \rtimes \hat{A}$ is topologically
trivial since we can construct cross-sections, but the cocycle is not
necessarily a coboundary.

\noindent{\bf 5.8 Example}. For the first description of the
Heisenberg group $G$ in Example 3.5 as an extension of ${\mathbb Z}^2$
by its center ${\mathbb Z}$, the action of ${\mathbb Z}^2$ on ${\mathbb T}\cong\hat{\mathbb Z}$ is trivial, therefore the groupoid ${\mathbb Z}^2\ltimes{\mathbb T}$ is just the cartesian product ${\mathbb Z}^2\times{\mathbb T}$ with $(a,b,z)$ composable with $(a',b',z')$ iff $z=z'$, and $(a,b,z)\cdot(a',b',z)=(a+a',b+b',z)$.
The cocycle $\omega:({\mathbb Z}^2\ltimes{\mathbb T})^2\to{\mathbb T}$ is given by \[\omega((a,b,z),(a',b',z))=z^{ab'}.\]
 Indeed, from the group operation in $G$, it follows that the map $f:{\mathbb Z}^2\times{\mathbb Z}^2\to{\mathbb Z}$ associated to the group extension is given by $f((a,b),(a',b'))=ab'$, and since ${\mathbb Z}\cong\hat{\mathbb T}$, we get the formula for $\omega$.

For the second description of $G$ as a semidirect product ${\mathbb Z}^2\rtimes{\mathbb Z}$, the cocycle is trivial since the extension splits, but the action of ${\mathbb Z}$ on $\widehat{\mathbb Z}^2\cong{\mathbb T}^2$ is given by
$\alpha(x,y)=(x,xy)$, and is induced by the conjugation with $[0,1,0]$ on ${\mathbb Z}^2$ and dualization.



\end{document}